\documentclass[12pt]{amsart}
\usepackage{amssymb, epsfig, latexsym}

\title{Polyomino Convolutions and Tiling Problems}
\author{Ali Ulas Ozgur Kisisel}

\newtheorem{Def}{Definition}
\newtheorem{Thm}{Theorem}
\newtheorem{Prop}{Proposition}

\begin{document}

\begin{abstract}
We define a convolution operation on the set of polyominoes and
use it to obtain a criterion for a given polyomino not to tile the
plane (rotations and translations allowed). We apply the criterion
to several families of polyominoes, and show that the criterion
detects some cases that are not detectable by generalized coloring
arguments.
\end{abstract}

\maketitle

\section{Introduction:}

Tiling properties of polyominoes have been studied by many authors
using various methods (see for instance \cite{Bru}, \cite{CL},
\cite{Go3}, \cite{Wal}).
We will assume that a polyomino $f$ is a map from $\mathbf{Z}\times \mathbf{Z%
}$ to $\mathbf{Z}$ that takes the values 0 and 1. We associate a value 1 of $%
\ f$ at $(n,m)\;$with the unit square $[n,n+1]\times \lbrack
m,m+1]$. This allows us to envision polyominoes in the usual way
as tiles in the plane. We assume $f$ \ has finitely many occupied
squares. (Several authors including S. Golomb refer to these
objects as quasi-polyominoes (\cite{Go2}, pg. 85), and reserve the
term polyomino for rookwise connected figures. A figure is
rookwise connected if it can be constructed by placing squares in
a way that each square except for the first shares an edge with a
previously placed square. We will not assume rookwise
connectedness or topological connectedness unless declared.)

\begin{Def}
Let $f$ be any map $\mathbf{Z}\times \mathbf{Z}\rightarrow
\mathbf{Z}$ \ (not necessarily a polyomino map). Suppose $f(n,m)$
is nonzero for finitely many pairs $(n,m)$. Let
$|f|_{1}=\sum_{n,m}|f(n,m)|$, and let $|f|_{\infty }$ be the
number of pairs $(n,m)$ such that $f(n,m)$ is not zero. Both of
these are norms in the usual sense. The assumption implies that
both $|f|_{1}$ and $|f|_{\infty }$ are finite, and we say that $f$
has finite area. If $f$ is a polyomino, the two norms are equal,
and we use the notation $|f|$.
\end{Def}

\begin{Def}
For any map $f$ as above, of finite area, let diam$(f)$ denote the
maximum of the distances between pairs of points of support of $f$
under the taxicab metric. The taxicab distance between two points
is the minimum number of grid steps from one point to the other (A
grid step is a move from $(x,y)$ to $(u,v)$ where $|x-u|+|y-v|=1$).
\end{Def}


\begin{Def}
Suppose $f$ and $g$ are two maps $\mathbf{Z}\times \mathbf{Z}\rightarrow
\mathbf{Z}$ such that at least one of them has finite area. Define $h=f\star
g$ as:
\begin{equation}
h(n,m)=\sum_{k,l}f(k,l)g(n-k,m-l)
\end{equation}
\end{Def}

$h(n,m)$ counts the number of intersections of $f$ with an $(n,m)$
translate of the reflection of $g$ across the origin. It is clear
that the sum in the definition is finite for a fixed $(n,m)$, and
$h$ has finite area if both $f$ and $g$ do. We call $h$ the
convolution of $f$ and $g$. It is easy to see that diam$(h)\leq $diam$(f)+$%
diam$(g)$.

If $f$ and $g$ are polyominoes, then their convolution $h=f\star
g$ is not a polyomino in general since $h$ may assume values other
than $0$ or $1$. But one may obtain a polyomino from $h$ by
reducing each $h(n,m)$ to $0$ or $1$ depending on its congruence
class modulo 2. In this way, we obtain a convolution operation on
the set of polyominoes. This operation inherits the associativity
and bilinearity properties of the usual convolution, since
reduction modulo 2 commutes with the operations of addition and
multiplication that constitute the convolution. Denote the
composition of convolution and the reduction modulo 2 by $f\star
_{2}g$. Similarly, denote the composition of convolution and
reduction modulo $n$ for an arbitrary modulus $n$ by $f\star
_{n}g$.

Here is our main observation:

\begin{Thm}
\label{thm:main} Suppose that $f$ is a polyomino symmetric under
rotations of $90$ degrees. Suppose $g$ is a polyomino. Then if
$|f\star _{n}g|_{1}<\bar{|f|}|g|$, or if $|f\star _{n}g|_{\infty
}<sgn(\bar{|f|})|g|$, where $\bar{|f|}$ denotes the unique integer
among $0,1,...,n-1$ congruent to $|f|$ modulo $n$, then copies of
$g$ cannot tile the plane (i.e. cover it without overlaps)\ ,
translations and rotations being allowed.(Here we are thinking of
$g$'s as tiles)
\end{Thm}

\noindent \textbf{Proof:} Say that copies of $g$ tile the plane.
This is another way to say that the full plane is a sum of non
overlapping translates of copies of $g$ and its rotations. Since
$f$ is rotationally symmetric, $f$ convolved with $g$ has the same
norms as $f$ convolved with a rotation of $g$. Consider a minimal
pattern of $g$'s in this tiling that contains a full $N$ by $N$
square. Call this figure $G$ for reference. $G$ is certainly
contained in a $N+2$diam$(g)$ by $N+2$diam$(g)$ square, otherwise
it wouldn't be minimal. Since $G$ is obtained as a disjoint sum of
$g$'s, its norm is simply the sum of the norms of its
constituents. We are going to estimate norms of $f\star _{n}G$
from two directions.

First of all, $f\star G$ has at least $(N-2$diam$(f))^{2}$ points of value
$|f|$, since at least that many translates of $f$ fall completely into the $N
$ by $N$ square. Reducing modulo $n$, we obtain $|f\star _{n}G|_{1}\geq (N-2$%
diam$(f))^{2}\bar{|f|}$ and $|f\star _{n}G|_{\infty }\geq (N-2$diam$(f))^{2}$%
sgn$(\bar{|f|})$. On the other hand $G$ is made up of at most $\frac{(N+2%
\mathrm{diam}(g))^{2}}{|g|}$ copies of $g$. If we assume to the
contrary that the inequalities in the hypothesis may hold, by the
 triangle inequality we obtain that
\begin{equation}
(N+2\mathrm{diam}(g))^{2}\frac{\bar{|f|}|g|-1}{|g|}\geq |f\star
_{n}G|_{1}\geq (N-2\mathrm{diam}(f))^{2}\bar{|f|}
\end{equation}
or
\begin{equation}
(N+2\mathrm{diam}(g))^{2}\frac{\mathrm{sgn}(\bar{|f|})|g|-1}{|g|}\geq
|f\star _{n}G|_{\infty }\geq (N-2\mathrm{diam}(f))^{2}\mathrm{sgn}(\bar{|f|})
\end{equation}
Both inequalities fail to hold asymptotically for large values of
$N$, since the coefficients of $N^{2}$ on the left hand sides of
the equations are strictly less than those on the right hand
sides. This contradiction finishes the proof. $\Box $

We remark that the theorem remains valid if we replace $g$ by a
finite collection of prototiles $g_{1},...,g_{k}$ such that the
inequalities hold for each of them separately.

\section{Some Applications}

We would like to show some applications of the criterion. Our first example
is a certain sequence of disconnected polyominoes.

\begin{figure}
\begin{center}
\epsfig{file=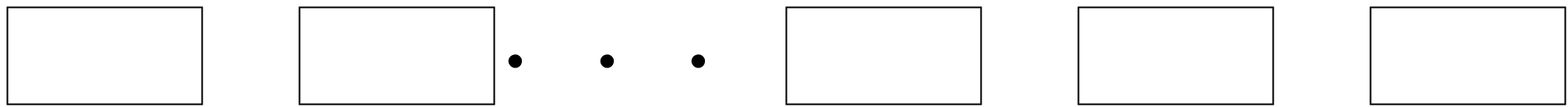,height=8cm,width=.5cm,angle=90} \caption{
$D_{n}$ } \label{fig:dn}
\end{center}
\end{figure}

We define a sequence $D_{n}$ of polyominoes as follows: $D_{n}$ is obtained
by aligning $n$ dominoes horizontally along their longer sides, leaving a
spacing of one square between any two consecutive dominoes (see figure \ref
{fig:dn}). For instance, an accordingly positioned $D_{n}$ would occupy the
squares $(0,0),(1,0),(3,0),(4,0),(6,0),(7,0),...,(3n-3,0),(3n-2,0)$. $D_{1}$
is a domino itself. Therefore it tiles the plane in many ways. The question
for $n\geq 2$ has the following answer:

\begin{Prop}
$D_{n}$ tiles the plane iff $n\leq 3$, translations and rotations allowed.
\end{Prop}

\begin{figure}
\begin{center}
\epsfig{file=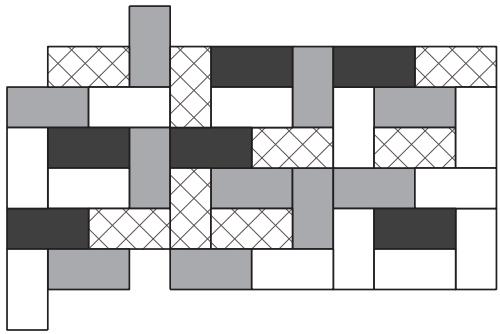,height=4.4cm,width=6.6cm,angle=0}
\caption{Tiling a $12\times 6$ torus by $D_{2}$'s} \label{fig:d2}
\end{center}
\end{figure}

\begin{figure}
\begin{center}
\epsfig{file=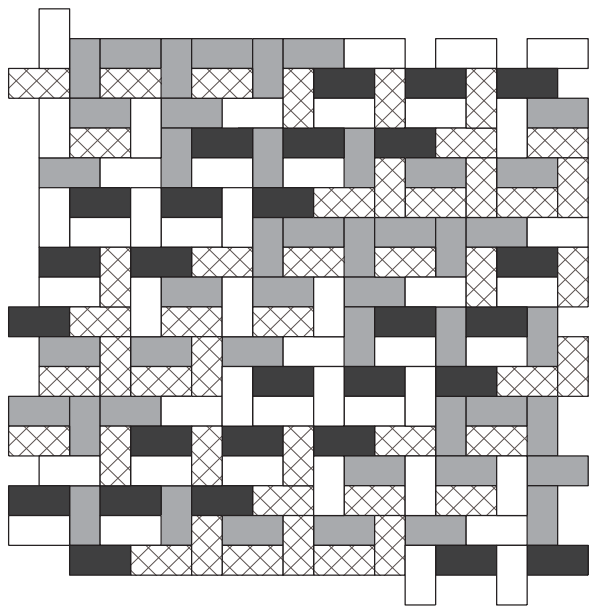,height=8cm,width=8cm,angle=0}
\caption{Tiling an $18\times 18$ torus by $D_{3}$'s}
\label{fig:d3}
\end{center}
\end{figure}

\noindent \textbf{Proof:} Examples of tilings for $D_{2}$ and $D_{3}$ are
shown in figures \ref{fig:d2} and \ref{fig:d3} respectively. Both tilings
are doubly periodic, thus only one fundamental region is shown in either
case. We must remark that a tiling pattern for $D_{2}$ or $D_{3}$ needs to
obey severe restrictions, and our guess is that the $D_{3}$ tiling is
essentially unique.

\begin{figure}
\begin{center}
\epsfig{file=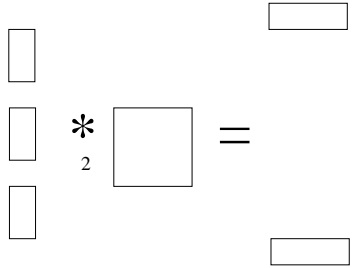,height=4.55cm,width=3.5cm,angle=270}
\caption{$D_{3}\star s_{3\times3}$} \label{fig:sd3}
\end{center}
\end{figure}

Next we prove the impossibility part of the assertion. Suppose that $%
S_{3\times 3}$ represents the $3$ by $3$ square polyomino. It is not hard to
check that  $|S_{3\times 3}\star _{2}D_{n}|=6$ for any value of $n$ (see
figure \ref{fig:sd3}). This happens since all but 6 translates of $%
S_{3\times 3}$ meet $D_{n}$ in an even number of squares. The 6 are those
where $S_{3\times 3}$ meets the first or last square of $\ D_{n}$.
Therefore, by Theorem \ref{thm:main} , $D_{n}$ cannot tile the plane if $%
|D_{n}|=2n>6$, i.e. if $n>3$. $\Box $

There are many ways that one can seek generalizations of this
example. A similar argument works for the negative part of the
corresponding assertion on higher dimensional analogues. We show
another generalization since it uses the other norm $|f|_{_{\infty
}}$:

\begin{Prop}
Let $D_{n,a,b}$ represent the polyomino obtained by aligning $n$ horizontal
bars of length $a$, leaving a spacing of $b$ blank squares between any two
consecutive bars (Therefore, the $D_{n}$ above are $D_{n,2,1}$ with this
notation). Then, if $b^{2}$ is not divisible by $a$, $D_{n,a,b}$ cannot tile
the plane if $n>\frac{2(a+b)(a-1)}{a}$.
\end{Prop}

\noindent \textbf{Proof:} Let $S_{(a+b)\times (a+b)}$ represent the square
polyomino of side length $(a+b)$. Then $|S_{(a+b)\times (a+b)}\star
_{a}D_{n,a,b}|_{\infty }=2(a-1)(a+b)$, and $|D_{n,a,b}|=na$. Thus the
inequality follows from the theorem unless $(a+b)^{2}$ is $0$ modulo $a$.
This is equivalent to $a|b^{2}$. $\Box $

\begin{figure}
\begin{center}
\epsfig{file=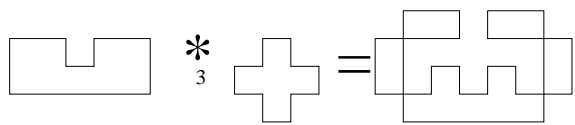,height=7.5cm,width=1.5cm,angle=270}
  \caption{A 9-omino which doesn't tile}
\label{fig:tile1}
\end{center}
\end{figure}

\begin{figure}
\begin{center}
\epsfig{file=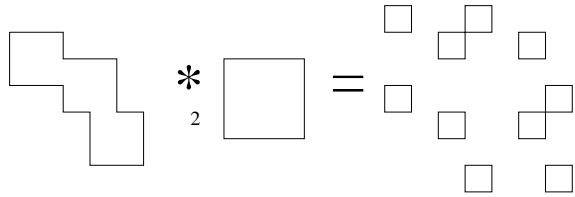,height=7.5cm,width=2.5cm,angle=270}
\caption{A 12-omino which doesn't tile} \label{fig:tile4}
\end{center}
\end{figure}

Next, we consider some rookwise connected polyominoes. All such
polyominoes of area 6 or less tile the plane (\cite{KV}), so we
have to consider larger polyominoes. The first polyomino in figure
\ref{fig:tile1} is a 9-omino that clearly doesn't tile. This is
provable by our criterion as well, as demonstrated in the same
figure. Figure \ref{fig:tile4} shows another non-tiler , and this
is also easy to prove directly.

We will call a polyomino $L$ a ``log'' if it is an $a$ by $b$
rectangle with $a>1$ and $b>1$. Let $L^{'}$ be the $a+2$ by $b+2$
rectangle containing $L$ in the middle. We define a ``log with
spikes'' (or a ``spiky log'') to be a polyomino obtainable from
such an $L$ by adjoining a number of 1 by 1 squares directly to
$L$ (each sharing an edge with a square of $L$)  so that there are
at least two blank squares between any two of them (Around
corners, count along the squares of $L^{'}-L$). We call the 1 by 1
squares ``spikes''.

\begin{Prop} No log with more than four spikes (moreover, no finite collection of
such prototiles) can tile the plane.
\end{Prop}

\begin{figure}
\begin{center}
\epsfig{file=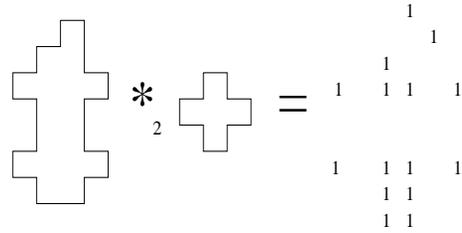,height=6cm,width=3cm,angle=270} \caption{A
spiky log with 5 spikes can't tile} \label{fig:tile3}
\end{center}
\end{figure}

\noindent \textbf{Proof:} Convolve the polyomino with the X
pentomino modulo 2. One may verify that convolving an $a$ by $b$
rectangle with the X pentomino gives a polyomino of norm $ab+8$,
and each spike placed on the log reduces the norm of the result by
1, while increasing the norm of the initial object by 1.
Therefore, if 5 or more spikes are placed, the convolution is norm
decreasing. (See the example in figure \ref{fig:tile3})$\Box$

We call a polyomino a ``snake'' if no subset of its squares is a T
tetromino or a square tetromino. We say that the snake makes $n$
``U-turns'' if it has $n$ distinct subsets forming U-pentominoes.

\begin{Prop} No snake making 3 or more U-turns (moreover, no finite collection of such prototiles)
can tile the plane.
\end{Prop}

\noindent\textbf{Proof:} Convolve the polyomino with the X
pentomino modulo 3 and look at the $|\phantom{x}|_{1}$ norm.
Except for the two squares at the ends, each square of the snake
has three neighbors (counting the square itself), so these do not
contribute. The squares at the ends may contribute 4 in total at
most. Any  square not on the snake makes a contribution only if it
shares an edge with the snake. If the snake has $n$ squares, the
maximum possible total contribution of this type is the number of
edges, $2n+2$. Every U-turn costs 3, therefore if there are more
than 2 U-turns, the norm of the convolution is less than $2n$, and
the criterion gives the result. $\Box$

 Golomb defines a
``reptile'' to be a polyomino that tiles a larger copy of itself
\cite{Go2}. All reptiles tile the plane. Therefore snakes making 3
U-turns (actually even 2 U-turns) cannot be reptiles!

\section{A Comparison to Coloring Arguments}
\
There are several other sufficient criteria to prove the
impossibility of tiling a given figure by another. Of these,
perhaps the best known are coloring and generalized coloring
arguments. One may ask where our criterion stands. We show that
there exist tiling problems such that the impossibility is
detected by our criterion whereas no generalized coloring argument
can do so. We follow the method in \cite{CL}. If a generalized
coloring argument proves impossibility of tiling $R$ with $f$,
then it also proves impossibility of a ``signed tiling'' of $R$
with $f$. A signed tiling permits using the map $-f$ as well as
$f$, and of course, overlaps allowed.

\begin{Thm}
(i) It is not possible to tile a torus by $D_{4}$'s.

(ii) There exist signed tilings of a $24$ by $12$ torus by $D_{4}$'s.
\end{Thm}

\noindent \textbf{Proof:} (i) is clear. If $D_{4}$ tiled a torus, it would
tile the plane in a doubly periodic way. But this was shown not to happen.

(ii) Notice that $f$ superposed with $-f$ shifted $3$ squares to the right
gives a map $g$ such that $g(0,0)=g(1,0)=1,g(12,0)=g(13,0)=-1$, and $0$
otherwise. Horizontally, stack $6$ $g$'s, with a shift of two squares
between consecutive $g$'s. We obtain a new map $h$ such that $h(k,0)=1$ for $%
k=0,...,11$, $h(k,0)=-1$ for $k=12,...,23$, and $0$ otherwise. Rotate $h$ $90
$ degrees clockwise.

Next, stack $12$ $D_{4}$'s vertically. We get a figure of $4$ rectangles of
dimensions $12$ by $2$, longer sides vertical, with one horizontal
separation between neighboring rectangles. Using copies of $h$, we can shift
the $2$nd and $4$th rectangles up by $12$ spaces while leaving the other two
untouched. A horizontal stack of three copies of this final figure gives a
torus tiling. $\Box$



\section*{Acknowledgement}

I wish to thank Izzet Pembeci, Ali Tamur and Ertem Tuncel for many
helpful conversations. The motivation for the main idea came from
their minesweeper puzzle with clues modulo 2. I also wish to thank
Terry Tao for helpful conversations, and Jeffrey Lagarias for
numerous useful suggestions and comments.

\begin{figure}
\begin{center}
\epsfig{file=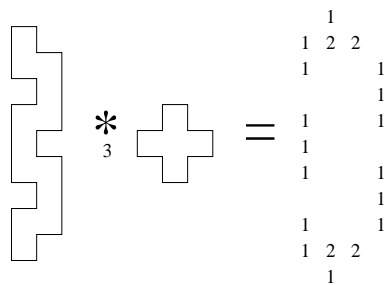,height=5cm,width=3.67cm,angle=270}
\caption{A snake making 3 U-turns can't tile} \label{fig:tile2}
\end{center}
\end{figure}

\end{document}